\documentclass[10pt,letterpaper]{article}
\usepackage[T1]{fontenc}
\usepackage[utf8]{inputenc}
\usepackage{geometry}
\usepackage{amsmath,mathtools,booktabs,tabularx,array,longtable}
\usepackage{caption,graphicx,needspace}
\usepackage[section]{placeins}
\usepackage[numbers]{natbib}
\usepackage{makecell}
\usepackage{enumitem}
\usepackage{multirow}
\usepackage{tikz}
\usepackage{colortbl}
\usetikzlibrary{arrows.meta,positioning}
\usepackage{authblk}
\usepackage{fontawesome5}
\usepackage{xurl}
\usepackage[pagebackref,hidelinks]{hyperref}
\hypersetup{
 pdftitle={How Much Can Large Language Models Help Solve MIPs?},
 pdfauthor={Yicheng Huang, Wenzhi Gao, Dongdong Ge, Madeleine Udell, Yinyu Ye},
 pdfsubject={LLM-assisted methods for mixed-integer programming}}
\newcommand*{\websitefont}{\fontfamily{FiraMono-TLF}\selectfont}

\newcolumntype{Y}{>{\raggedright\arraybackslash}X}
\newcommand{\inst}[1]{\texttt{\detokenize{#1}}}
\newcommand*{\miprankingentry}[1]{\textsf{#1}}
\newcommand{\llm}{\href{https://llm4mip.github.io}{{\websitefont LLM4MIP}}}

\begin{document}

\title{How Much Can Large Language Models Help Solve MIPs?}
\author[1,*]{Yicheng Huang\footnote{huangyc.opt@gmail.com. The first two authors co-lead the project.}}
\author[2,*]{Wenzhi Gao\footnote{gwz@stanford.edu}}
\author[3]{Dongdong Ge\footnote{dongdong@shanshu.ai}}
\author[2]{Madeleine Udell\footnote{udell@stanford.edu}}
\author[2]{Yinyu Ye\footnote{yyye@stanford.edu}}
\affil[1]{Shanghai University of Finance and Economics}
\affil[2]{Stanford University}
\affil[3]{Cardinal Operations}
\date{November 27, 2026\\[3pt]{\normalsize Initial version (\texttt{v0}): October 7, 2026}\\[7pt]
  {\large
  \href{https://github.com/LLM4MIP/llm4mip.github.io}{\faGithub\enspace GitHub}\qquad
  \href{https://llm4mip.github.io/}{\faGlobe\enspace Website}}}
\setcounter{footnote}{1}
\maketitle

\begingroup
\renewcommand{\thefootnote}{\fnsymbol{footnote}}

\begin{abstract}
Large language models (LLMs) have recently made dramatic progress on problems in theoretical mathematics, including many that were previously open. Can they also help solve open instances of mixed-integer programming (MIP)? Solving a difficult MIP typically requires a combination of mathematical reasoning, numerical computation, and optimization software. Although LLMs have been studied for a range of operations research tasks, relatively little work has directly evaluated their ability to solve MIP instances. 
In this technical report, we introduce {\llm}, a systematic benchmark for existing LLMs on 217 open instances from \texttt{MIPLIB} 2017.
Our results show that LLMs can accelerate the solution process by identifying and exploiting instance-specific structure. They use this structure to design heuristics, derive cutting planes, and even sometimes construct mathematical proofs. These results suggest that LLMs have the potential to significantly lower the barrier to developing customized solution methods for difficult MIPs.
\end{abstract}

\begingroup
\setlength{\parskip}{3pt plus 1pt}
\section{Introduction}
Mixed-integer programming {\cite{wolsey2020integer,sahinidis2019mixed}}
is arguably the most widely used mathematical optimization framework, with
applications throughout production planning {\cite{pochet2006production}},
combinatorial optimization {\cite{wolsey1999integer}}, energy management
{\cite{wood2013power}}, and revenue management {\cite{benati2007mixed}}.
Although MIP is a powerful modeling tool, MIP instances are often hard to solve even
for state-of-the-art commercial solvers
{\cite{copt,gurobi,cplex2009v12,optimizer2021v8}}. To solve an instance to
optimality, a solver must find a feasible solution and certify that no better
one exists, and a general-purpose solver implements a number of heuristics and
cutting planes to adapt to instance-specific structures. Once a new structure
is understood, it can be internalized by the solver and generalized to future
instances upon detection. Modern commercial solvers have built up a pool of
such structures through years of exposure to diverse customer instances.
Despite years of development, solvers can still struggle when facing new MIP
instances with unseen structures. In the pre-LLM era, structure investigation
required multiple areas of expertise, including mathematical reasoning, algorithm
design, and familiarity with optimization software. Consequently, even for
a single instance, identifying a useful structure and turning it into a method
that speeds up the solve can take considerable human effort.

Recently, large language models have demonstrated striking capabilities
in mathematical reasoning. Many open problems across theoretical and applied
mathematics have seen progress with the help of LLMs
{\cite{openai2026unitdistance,openai2026tenadvances,anthropic2026blackholes,alpoge2026riemann,anthropic2026riemannprogress}}.
On the other hand, LLMs have also long been recognized as powerful tools for
coding and algorithm design
{\cite{chen2021evaluating,romera2024mathematical,alphaevolve}}. Given that
some recent studies have also demonstrated LLMs' potential on operations
research tasks {\cite{li2026constructing,kong2026frontieror}}, modern LLMs
appear to have most of the expertise needed to identify useful structure in
MIPs. Given the limited research directly benchmarking LLMs on solving MIP
instances, a natural question arises:

\begin{center}
\textit{{Can LLMs help solve MIPs?}}	
\end{center}

We answer this question in the affirmative through the {\llm} project.

\begin{itemize}[leftmargin=10pt,topsep=0pt,itemsep=0pt]
  \item We conduct a systematic benchmark of LLMs' ability to
  solve real-world MIPs, covering 217 open \texttt{MIPLIB} instances \cite{gleixner2021miplib} (10/7/2026). Our results show that LLMs can identify useful
  instance-specific structure and exploit it to design heuristics, derive cutting planes, and even develop proofs tailored to individual instances.

  \item We distill reusable skills from our work on the 217 \texttt{MIPLIB}
  instances and openly release artifacts (\url{https://llm4mip.github.io}) produced during the
  LLM-assisted solution process. A comparison on the test set
  shows that these skills can transfer to other MIP instances.
\end{itemize}

\subsection{Related work}

There is a growing body of research on LLMs for operations research and optimization. Due to space limits, we focus on closely related work
and refer interested readers to more comprehensive surveys
{\cite{xiao2025survey,guo2026natural,tian2026advances,liu2026systematic,da2026large,xiu2026large}}.

\paragraph{LLMs for operations research.}

Much recent work on LLMs for OR
focuses on optimization modeling, with agent-based frameworks
{\cite{ahmaditeshnizi2024optimus,xiao2024chain}}, fine-tuning approaches
{\cite{huang2025orlm}}, and reinforcement learning approaches
{\cite{shu2025llmopt}}. Other works use LLMs to propose candidate solutions or design heuristics
for combinatorial optimization {\cite{yang2024opro,eoh,reevo}}. Another notable recent example is FrontierOR, which evaluates the solution quality and runtime of LLM-designed algorithms
{\cite{kong2026frontieror,wang2026frontieror}}.

\paragraph{LLMs for optimization solvers.} A particularly related line of
research investigates how LLMs can assist existing solver-based workflows.
Machine learning has been used to guide branching
{\cite{khalil2016learning,gasse}}, find primal solutions
{\cite{nair2020solving,zhang2026shsp}}, and select cutting planes
{\cite{paulus2022learning}}. Recent work uses LLMs to design solver
procedures such as branching policies
{\cite{hou2026llmbranch,zhang2026bife,zhou2026large}}. LLMs have also
been used to configure cutting-plane separators {\cite{lawless2025llms}}
and to tune parameters for individual instances, often using solver logs
{\cite{you2026log,luo2026beyond,luo2026grimip,zhao2026nominate}}. Other
applications include MIP re-optimization and SMT solver tuning
{\cite{liuonline,wang2026agent}}. Some works use LLMs to generate families of
inequalities or to search for solver-efficient MIP formulations
{\cite{yazdani2025evocut,formuevo}}. A few recent projects also
use LLMs to write numerical optimization solvers from scratch
{\cite{korda2026brisk,polkowski_popso}}, each achieving promising performance compared to open-source solvers.
 
\paragraph{Concurrent results.} After our project website was online, Anthropic applied similar ideas to open \texttt{MIPLIB} instances with Claude agents and reported
improvements to best-known feasible solutions {\cite{claude2026miplib}}.
\texttt{Gurobi} also benchmarked this approach on \texttt{MIPLIB} instances and provided insights on how LLMs can help MIP solving \cite{rothberg2026llmsmip}.

\paragraph{Structure of the report.} The rest of this technical report is organized as follows. In \textbf{Section \ref{label:method}},  we discuss our prompting approaches and the skill-distillation process. \textbf{Section \ref{sec:exp}} covers the detailed experiment results and representative case studies. We discuss the potential impact of LLMs on future MIP solvers in \textbf{Section \ref{sec:discussion}}.

\section{Large language model-assisted MIP solving} \label{label:method}

This section introduces the LLM-based workflow that is used throughout our
experiments.

\subsection{Initial screening and preprocessing}

We start by extracting the detailed statistics for 217 open \texttt{MIPLIB} instances
from \url{https://miplib.zib.de/tag\_open.html} in \texttt{CSV} format.
Then we use a single prompt to let the LLM go through these instance statistics
and determine
\begin{itemize}
  \item a ranking of open instances by the confidence (out of 100) for solving
  each of them, and 
  
  \item the possible entry point of the improvement (primal or dual).
\end{itemize}
The LLM is allowed to access the \texttt{MIPLIB} website during the initial screening
process. After the initial screening, 17 instances are considered to have a high
confidence score ($\ge 60$), 57 to have a medium score (45--59), and the remaining
143 to have a low score ($<45$). \textbf{Figure~\ref{fig:screening-scores}} shows the
distribution of the confidence scores.
\textbf{Appendix~\ref{app:screening-scores}} lists the scores for individual
instances. The screening procedure is limited to 3600 seconds.

\begin{figure}[htbp]
  \centering
  \includegraphics[scale=1.0]{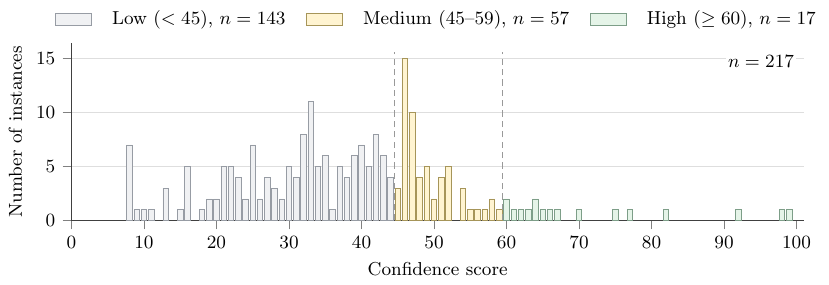}
  \caption{Distribution of heuristic confidence scores for the 217 open \texttt{MIPLIB}
  2017 instances after initial screening. Each bar shows the number of
  instances at one integer score.}
  \label{fig:screening-scores}
\end{figure}

\subsection{Per-instance prompting and tools}

Given a MIP instance in \texttt{MPS} format {\cite{gleixner2021miplib}},
detailed instance statistics are collected from \texttt{MIPLIB} (\textbf{Figure~\ref{fig:allcolor58}}). After that, the LLM is asked
to solve the instance based on all available information and tools, including

\begin{itemize}[label={},leftmargin=1em,itemsep=2pt,topsep=3pt]
  \item \textit{Results from previous screening.} The LLM is given the available
  best-known solution from \texttt{MIPLIB}. It also has access to the confidence
  score and the entry point information obtained from the screening process.

  \item \textit{Internet and literature search.} The LLM is allowed to search online for
  additional details on the problem background and previous attempts, including
  web forum discussions and related research papers.

  \item \textit{Numerical solvers.} The LLM is allowed to invoke optimization solvers
  including \texttt{Gurobi} and \texttt{COPT}, and also has the freedom
  to install additional open-source solvers (e.g., \texttt{OR-Tools}) in the
  workspace.
\end{itemize}

\begin{figure}[ht]
  \centering
  \includegraphics[height=.216\linewidth]{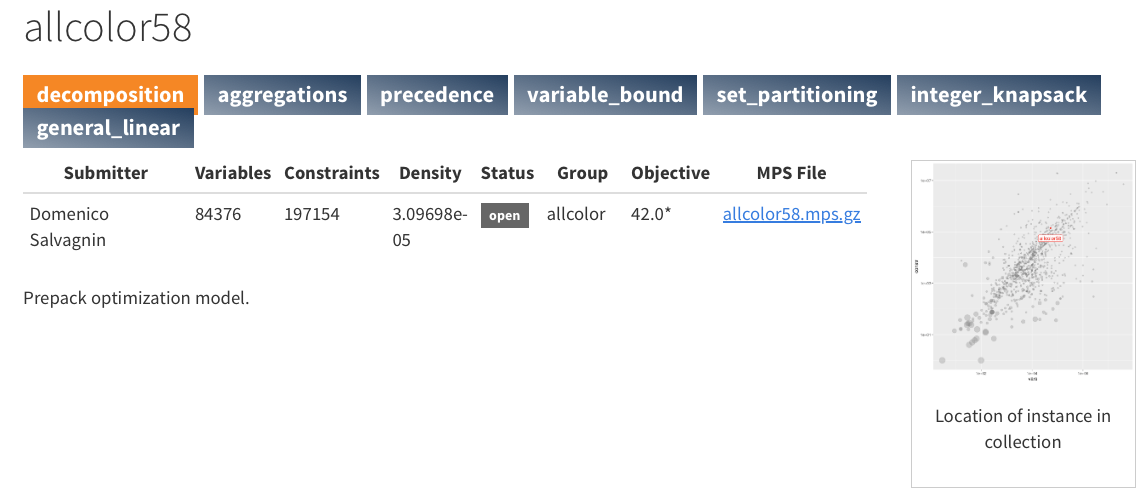}\hspace{.025\linewidth}\includegraphics[height=.216\linewidth,trim=8bp 0bp 0bp 7bp,clip]{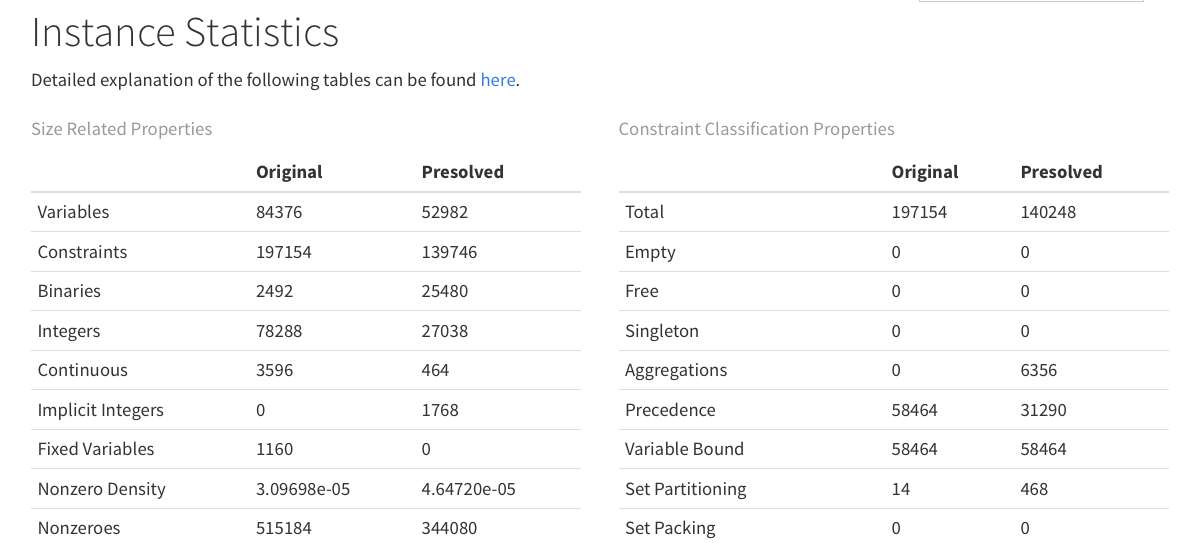}
  \caption{\texttt{MIPLIB} instance statistics for \texttt{allcolor58}.\label{fig:allcolor58}}
\end{figure}

The prompt is simple:

\smallskip
\begingroup
\setlength{\fboxsep}{10pt}
\setlength{\fboxrule}{.45pt}
\noindent\fcolorbox{black!22}{black!2}{\begin{minipage}{\dimexpr\linewidth-2\fboxsep-2\fboxrule\relax}
    \renewcommand*{\ttdefault}{FiraMono-TLF}
    \small\ttfamily\frenchspacing\raggedright
    \setlength{\parskip}{5pt}
    Please solve the open \texttt{MIPLIB} instance \textcolor{blue!40!black}{\{NAME\}}.

    The previous screening result is given in screening.csv.

    Feel free to look for additional information online and use solvers.
    \textcolor{blue!40!black}{\{SOLVER\}} is installed locally. Try exploiting
    the instance structure and use solvers sparingly.
  \end{minipage}}
\par
\endgroup
\smallskip

To prevent an LLM from calling a solver directly without exploiting structure, we explicitly ask it to use solvers sparingly. In the standard setup, each LLM-based solve is allowed a time limit of 7200 seconds.

\subsection{Result validation} \label{sec:validation}

\begingroup
\setlength{\parskip}{6pt plus 1pt}

After attempts are made on an instance, there are the following possible cases of
claimed improvements. The first two kinds of improvements are made on the
upper bound, and they can be mechanically verified.

\textbf{Primal bound improvement.}
Validating a primal solution can be done by plugging in the claimed solution
and verifying that no constraint is violated beyond the reported tolerances. We use \texttt{Gurobi} and
\texttt{COPT} for validation: a solution is considered valid when
both solvers can accept the solution as a valid warm start. Then the new primal
solution is recorded.

\textbf{Unbounded instance.}\ This case is rare (only on
\texttt{supportcase35}), and the LLM provides an improving ray for the problem. Concretely, given an instance with objective vector $c$, the LLM
finds a direction $d$ and a feasible solution $x$ such that $x + \alpha d$
is feasible for all $\alpha \geq 0$ and $\langle c, d \rangle < 0$.
Then the objective tends to
$-\infty$ as $\alpha \to +\infty$.

The last two types of improvements are made on the lower bound, and they are
much more challenging to verify.

\textbf{Dual bound improvement.}\ LLMs can report the dual bound progress arising from different
sources, including reformulation, customized cutting plane/valid constraint
generation, logical reasoning, and theorem proof.

The trustworthy case is when a dual improvement can be certified by a
solver running on the original instance. Unfortunately, it is almost never the
case, and most of the dual bound claims require human verification.

To reduce the effort of manual verification, we
\begin{itemize}[leftmargin=14pt,topsep=3pt,itemsep=2pt]
  \item cross-check the claimed dual bound improvement independently with
  different LLMs.

  \item ask the LLM to generate a \texttt{Lean} proof for the proposed cutting
  planes/theorems.
\end{itemize}
Humans then manually verify these results. See the
experiment section for the progress of manual verification.

\textbf{Infeasible instance.}\ Infeasible instances have similar
verification barriers to dual bound improvements.
\par
\endgroup

\subsection{(Iterative) skill distillation}

After solving the instances, we save the detailed reasoning process of LLMs and distill structures that can be used on other instances. This distillation process yields skill packages that can be conveniently loaded into LLMs. We distill two skill sets for primal and dual bound improvements, respectively; each skill tells the LLM what structure to look for, how to exploit it, and how to verify the result on the original instance. 

\textbf{Iterative distillation.} Whenever a few more instances are solved, the skills can be iteratively updated. We mock the process of solver development by adding useful structures and recording failed attempts from the new solves. This iterative process makes it easier for the LLM to solve future instances with similar structures (\textbf{Figure~\ref{fig:skill-distillation}}). In our experiments, we split the 217 \texttt{MIPLIB} instances into a first group of 132 and a second group of 85. The group of 132 instances contains most instances with high-to-medium confidence scores. We then use the skills
from the first 132 instances to help solve the remaining 85, which have lower
confidence scores overall.

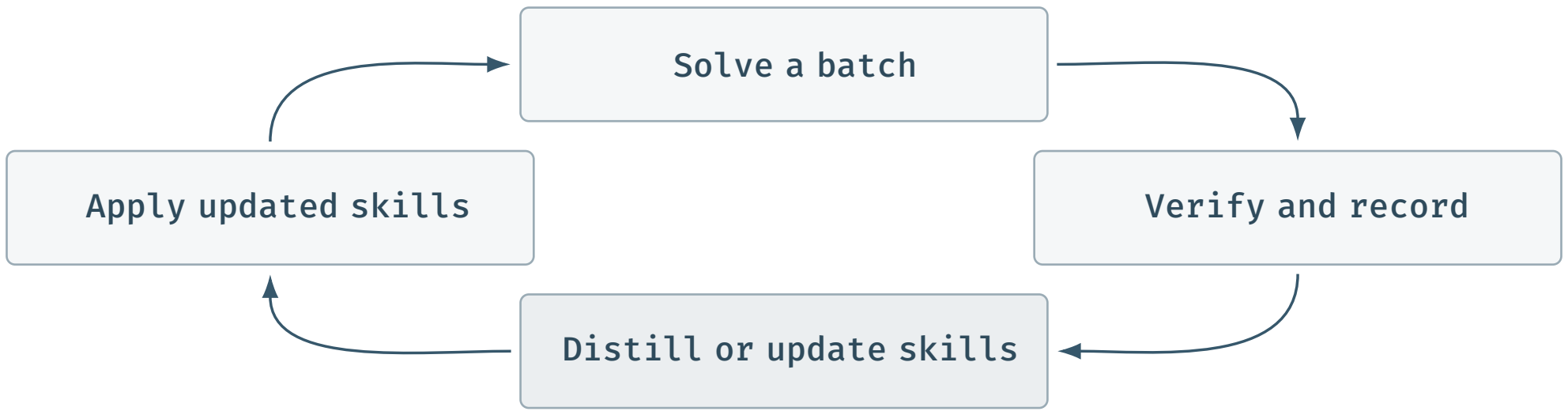
\begin{figure}[!htb]
  \centering
  \definecolor{skillworkflow}{HTML}{36576B}
  \begin{tikzpicture}[
    every node/.style={font=\footnotesize\fontfamily{FiraMono-TLF}
      \fontseries{medium}\selectfont,text=skillworkflow!85!black},
    stage/.style={draw=skillworkflow!50,fill=skillworkflow!5,rounded corners=2pt,
      line width=.5pt,text width=42mm,align=center,inner xsep=3pt,
      minimum height=9.5mm},
    flow/.style={-{Latex[length=2mm,width=1.4mm]},draw=skillworkflow,
      line width=.65pt,shorten <=2pt,shorten >=2pt},
  ]
    \node[stage] (solve) at (0,12mm) {Solve a batch};
    \node[stage] (verify) at (43mm,0) {Verify and record};
    \node[stage,fill=skillworkflow!10] (update) at (0,-12mm)
      {Distill or update skills};
    \node[stage] (apply) at (-43mm,0) {Apply updated skills};

    \draw[flow] (solve.east) to[out=0,in=90] (verify.north);
    \draw[flow] (verify.south) to[out=-90,in=0] (update.east);
    \draw[flow] (update.west) to[out=180,in=-90] (apply.south);
    \draw[flow] (apply.north) to[out=90,in=180] (solve.west);
  \end{tikzpicture}
  \caption{Iterative skill distillation.}
  \label{fig:skill-distillation}
\end{figure}

\section{Experiments} \label{sec:exp}
\paragraph{Experimental setup.}
In the standard setup, we use \texttt{Codex GPT 5.6-sol} as the base language model and run a single \texttt{Codex} session per solve attempt. Each attempt has a wall-clock budget of 7200 seconds, including time for structural analysis, solver calls, and solution validation. MIP solver calls use \texttt{COPT} or \texttt{Gurobi} by default.

The primal baseline is the objective value of the best-known solution reported in \texttt{MIPLIB} v36 (see \textbf{Figure~\ref{fig:allcolor58}}). For dual improvements, the baseline is the final dual bound reported by \texttt{COPT} with a 10-hour time limit per instance. Following \textbf{Section~\ref{sec:validation}}, we certify optimality through mathematical proofs or solver-based numerical verification. In numerical verification, classification as optimal within tolerance requires a relative optimality gap below $10^{-6}$, where the gap between primal bound $P$ and dual bound $D$ is defined as

\begin{equation}
  \textsf{Gap}(P,D)\coloneqq\frac{P-D}{\max\{1,|P|,|D|\}}
  \label{eq:gap}
\end{equation}

for a minimization problem.

\subsection{Results on open \texttt{MIPLIB} instances}
\looseness=-1
\textbf{Table~\ref{tab:miplib}} summarizes the results on \texttt{MIPLIB}\footnote{The detailed instance results are available at \url{https://llm4mip.github.io/}.}. The overall performance of the workflow is surprising: LLMs achieve primal-bound improvements on 59 instances and dual-bound improvements on 136.
Together, these improvements resolved 49 open \texttt{MIPLIB} instances to optimality, with two additional infeasible instances and one unbounded instance. For the instances classified as optimal, 25 already have matching feasible objective values in \texttt{MIPLIB}. Therefore, their bottleneck is on the dual side.

\begin{table}[!htbp]
\centering\small
\caption{Results on 217 \texttt{MIPLIB} instances}
\label{tab:miplib}
\definecolor{outcomeoptimal}{HTML}{587F96}
\definecolor{outcomeinfeasible}{HTML}{75A98D}
\definecolor{outcomeunbounded}{HTML}{C49757}
\definecolor{outcomeunresolved}{HTML}{E3E8EC}
\newcommand*{\outcomekey}[1]{\textcolor{#1}{\rule{.7em}{.7em}}\enspace}
\begin{minipage}[c]{.37\linewidth}
  \small
  \renewcommand{\arraystretch}{1.0}
  \setlength{\aboverulesep}{.3ex}
  \setlength{\belowrulesep}{.45ex}
  \begin{tabular*}{\linewidth}{@{\extracolsep{\fill}}lrr@{}}
    \toprule
    Result & Count\\
    \midrule
    Instances studied & 217\\
    \addlinespace[1pt]
    \outcomekey{outcomeoptimal}Optimal & 49\\
    \outcomekey{outcomeinfeasible}Infeasible & 2\\
    \outcomekey{outcomeunbounded}Unbounded & 1\\
    \outcomekey{outcomeunresolved}Unresolved & 165\\
    \midrule
    New or improved solutions & 59\\
    Improved dual bounds & 136\\
    \bottomrule
  \end{tabular*}
\end{minipage}\hspace{.04\linewidth}\qquad\qquad
\begin{minipage}[c]{.31\linewidth}
  \centering
  \begin{tikzpicture}[
    callout/.style={-{Latex[length=1.2mm,width=.9mm]},line width=.45pt},
    callout label/.style={anchor=west,font=\scriptsize,text=black,inner xsep=2pt},
  ]
    \foreach \outcome/\outcomecount [remember=\sectorend as \sectorstart (initially 90)] in
      {outcomeoptimal/49,outcomeinfeasible/2,outcomeunbounded/1,outcomeunresolved/165}{
      \pgfmathsetmacro{\sectorend}{\sectorstart-360*(\outcomecount/217)}
      \fill[\outcome] (0,0) -- (\sectorstart:1.48cm)
        arc[start angle=\sectorstart,end angle=\sectorend,radius=1.48cm] -- cycle;
    }
    \fill[white] (0,0) circle[radius=.98cm];
    \node[align=center,inner sep=0pt] {\large 24.0\%\\[-1pt]\footnotesize 52 resolved};
    \pgfmathsetmacro{\infeasibleangle}{90-360*((49+2/2)/217)}
    \pgfmathsetmacro{\unboundedangle}{90-360*((49+2+1/2)/217)}
    \draw[callout,draw=outcomeinfeasible!70!black] (\infeasibleangle:1.48cm)
      -- (1.72cm,.43cm) -- (1.98cm,.43cm) node[callout label] {2 infeasible};
    \draw[callout,draw=outcomeunbounded!70!black] (\unboundedangle:1.48cm)
      -- (1.72cm,-.30cm) -- (1.98cm,-.30cm) node[callout label] {1 unbounded};
  \end{tikzpicture}
\end{minipage}
\end{table}

\paragraph{Sources of primal and dual improvements.}
\textbf{Table~\ref{tab:miplib-methods}} groups the improvements by the main method used to find the final feasible solution or establish the final dual bound. The detailed explanation for each solution method is summarized below. Perhaps the most surprising fact is that LLMs demonstrate strong capability of approaching MIP instances from very different perspectives, ranging from common techniques such as local search and warm-starting, to more complicated techniques such as reformulation and providing customized cutting planes. Some instances are solved with direct mathematical reasoning without relying on a solver. The experiments demonstrate LLMs' superior performance in exploiting instance-dependent structures.

\newcommand*{\instartifact}[1]{\href{https://llm4mip.github.io/instances/details/#1/#1-findings.tar.gz}{\inst{#1}}}
\begin{table}[!htbp]
\centering\small
\definecolor{methodbar}{HTML}{587F96}
\newcommand*{\methodlink}[2]{\begingroup
  \setlength{\fboxsep}{1pt}\setlength{\fboxrule}{.35pt}\hyperlink{method:#1}{\fbox{#2}}\endgroup
}
\newcommand*{\methodcount}[1]{\makebox[1.20cm][l]{\textcolor{methodbar}{\rule[\dimexpr .5ex-.21em\relax]{\dimexpr .38mm*#1\relax}{.42em}}}\makebox[.43cm][r]{#1}}
\newsavebox{\improvementheadingbox}
\newcommand*{\improvementheading}[1]{\sbox{\improvementheadingbox}{\textbf{#1}}\leavevmode\vtop{\hbox{\usebox{\improvementheadingbox}}\vskip\aboverulesep
    \hrule width\wd\improvementheadingbox height\lightrulewidth
    \vskip\belowrulesep
  }}
\captionsetup{skip=4pt}
\caption{Sources of improvement on \texttt{MIPLIB} instances}
\begin{tabularx}{\linewidth}{@{}>{\raggedright\arraybackslash}p{0.39\linewidth} c Y@{}}
\toprule
Method (click for explanation) & Instances & Representative instances (click \texttt{instance} to download)\\
\midrule
\multicolumn{3}{@{}l}{\improvementheading{Primal improvements (59 instances)}}\\
\methodlink{local-search}{Local search} & \methodcount{31} & \instartifact{nj1}, \instartifact{allcolor58}\\
\methodlink{primal-reformulation}{Reformulate or decompose the model} & \methodcount{13} & \instartifact{liu}, \instartifact{bharat}\\
\methodlink{solution-repair}{Repair a known solution} & \methodcount{8} & \instartifact{circ10-3}, \instartifact{fhnw-schedule-paira200}\\
\methodlink{logical-reasoning}{Direct logical reasoning} & \methodcount{3} & \instartifact{bmocbd3}, \instartifact{neos-3208254-reiu}\\
\methodlink{warmstarting}{Warm-starting a standard MIP solver} & \methodcount{2} & \instartifact{nag}, \instartifact{ns1905797}\\
\methodlink{literature-solution}{Look for a solution from the literature} & \methodcount{2} & \instartifact{fastxgemm-n3r21s3t6}, \instartifact{fastxgemm-n3r22s4t6}\\
\midrule
\multicolumn{3}{@{}l}{\improvementheading{Dual-bound improvements (136 instances)}}\\
\methodlink{relaxation}{Bounds from relaxations or subproblems} & \methodcount{31} & \instartifact{ns930473}, \instartifact{fhnw-schedule-paira200}\\
\methodlink{cuts}{Add valid inequalities (cuts)} & \methodcount{25} & \instartifact{fastxgemm-n3r21s3t6}, \instartifact{eva1aprime6x6opt}\\
\methodlink{counting}{Prove a bound by counting or arithmetic} & \methodcount{24} & \instartifact{liu}, \instartifact{sorrell7}\\
\methodlink{dual-reformulation}{Reformulate the model} & \methodcount{21} & \instartifact{momentum3}, \instartifact{cmflsp60-36-2-6}\\
\methodlink{objective-branching}{Branch on the objective value} & \methodcount{21} & \instartifact{graph40-40-1rand}, \instartifact{bppc6-06}\\
\methodlink{solver-bound}{Global bound from a standard solver} & \methodcount{12} & \instartifact{cmflsp40-36-2-10}, \instartifact{ger50-17-trans-pop-3t}\\
\methodlink{prior-bound}{Use a bound from prior work} & \methodcount{2} & \instartifact{circ10-3}, \instartifact{a2864-99blp}\\
\bottomrule
\end{tabularx}
\label{tab:miplib-methods}
\end{table}
\FloatBarrier

\begin{minipage}{\linewidth}
\textbf{Primal bound improvements.}
\begin{itemize}[label={},leftmargin=0pt,topsep=2pt,itemsep=2pt,parsep=0pt]
  \item \hypertarget{method:local-search}{\textit{Local search.}} Keep most binary variables fixed and try changes to a selected subset and solve a sub-MIP.
  \item \hypertarget{method:primal-reformulation}{\textit{Reformulate or decompose the model.}} Solve a more compact formulation of the original problem or a set of smaller problems when the problem can be decomposed.
  \item \hypertarget{method:solution-repair}{\textit{Repair a known solution.}} Extrapolate a known solution from a related model and repair violations. For example, the solution for \instartifact{fhnw-schedule-paira200} can be adapted to other instances from the \inst{fhnw-schedule} family.
  \item \hypertarget{method:logical-reasoning}{\textit{Direct logical reasoning.}} Reason about what the constraints describe and use those rules to put together a feasible solution. For \instartifact{neos-3208254-reiu}, a binary tree was built, and its parent and ancestor relationships were translated into values for the original variables.
  \item \hypertarget{method:warmstarting}{\textit{Warm-starting a standard MIP solver.}} Use the feasible solution found by the solver's own search. It can start from an existing candidate and continue looking for improvements.
  \item \hypertarget{method:literature-solution}{\textit{Look for a solution from the literature.}} Modify an existing construction from the literature to fit the instance.
\end{itemize}
\end{minipage}

\smallskip
\textbf{Dual bound improvements.}
\begin{itemize}[label={},leftmargin=0pt,topsep=2pt,itemsep=2pt,parsep=0pt]
  \item \hypertarget{method:relaxation}{\textit{Bounds from relaxations or subproblems.}} Solve a relaxation of the original problem or a few subproblems.
  \item \hypertarget{method:cuts}{\textit{Add valid inequalities (cuts).}} Add constraints that every feasible solution already satisfies.
  \item \hypertarget{method:counting}{\textit{Prove a bound by counting or arithmetic.}} Use the background of the original problem to reason how good any feasible solution can be. For \instartifact{liu}, the total area to be packed and the possible side lengths imply a lower bound on the size of the enclosing square.
  \item \hypertarget{method:dual-reformulation}{\textit{Reformulate the model.}} Rewrite the model and establish that its best objective value is preserved. A bound obtained from the new formulation then applies to the original problem.
  \item \hypertarget{method:objective-branching}{\textit{Branch on the objective value.}} Add a constraint of the form $\langle c, x \rangle \leq z$ and prove the resulting model infeasible.
  \item \hypertarget{method:solver-bound}{\textit{Global bound from a standard solver.}} Use the bound a MIP solver reports for the original instance.
  \item \hypertarget{method:prior-bound}{\textit{Use a bound from prior work.}} Match the instance to a problem already studied in the literature.
\end{itemize}

\paragraph{Can LLMs predict whether an instance can be closed?}
We compare the initial confidence scores with the final outcomes for the same 217 instances. As shown in \textbf{Figure~\ref{fig:confidence-resolution}}, the solve rate rises from 17.5\% (25 of 143) in the low-confidence group to 29.8\% (17 of 57) in the medium-confidence group and 58.8\% (10 of 17) in the high-confidence group. The low-confidence group accounts for nearly half of all resolutions (25 of 52) since it contains many more instances. Therefore, the experiment suggests that even with a few minutes of study, LLMs can already roughly judge the possibility of closing a hard MIP instance.

\begin{figure}[!htbp]
  \centering
  \includegraphics[width=\linewidth]{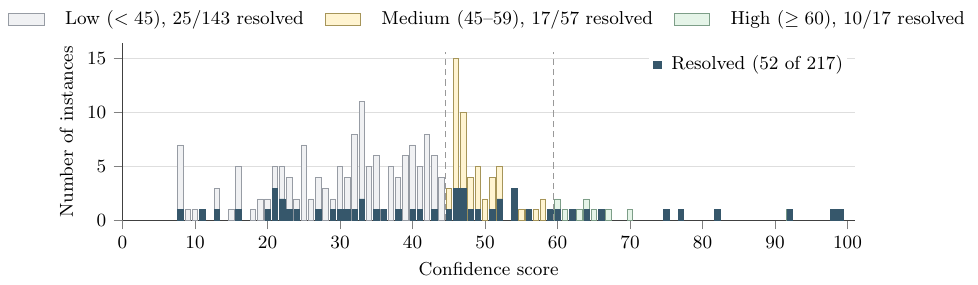}
  \caption{Total and resolved instances by initial confidence score. The dark bars show how many of them were resolved, and the dashed lines separate the confidence bands.}
  \label{fig:confidence-resolution}
\end{figure}
\FloatBarrier

\subsection{Distilled skill comparison}
After distilling the initial versions of the primal and dual skills, we evaluate them on 20 unseen open instances from the \texttt{MIPLIB} collection dataset.\footnote{Further details on the skill comparison are available at \url{https://llm4mip.github.io/skill/}.} Each skill run has a time limit of 9000 seconds.

\textbf{Table~\ref{tab:skills}} summarizes the comparison.
Compared with the LLM without skills, the primal skill improves the objective on 10 instances, ties on nine, and loses on one. The dual skill wins 17, ties 1, and loses 2. Compared to direct solving with \texttt{Gurobi}, the skills win on 17 of the instances.

\begin{table}[ht]
\centering
\small
\caption{Skill comparison on 20 instances.}
\setlength{\tabcolsep}{8pt}
\begin{tabular}{@{}l c c c c@{}}
\toprule
\multirow{2}{*}{Compared with}
& \multicolumn{3}{c}{Comparison results (wins/ties/losses)}
& \multirow{2}{*}{\shortstack{Mean gap\\reduction (pp)}} \\
\cmidrule(lr){2-4}
& Primal & Dual & Gap & \\
\midrule
LLM without skills
& 10/9/1 & 17/1/2 & 18/0/2 & 15.8079 \\

\texttt{Gurobi} (150 minutes)
& --- & --- & 17/1/2 & 19.2762 \\
\bottomrule
\end{tabular}
\label{tab:skills}
\end{table}
\FloatBarrier

Our results preliminarily show that the skills often generalize to new MIP instances.

\Needspace{18\baselineskip}
\subsection{Can LLMs replace MIP solvers?}
\label{sec:solver-replacement}

Interestingly, of the 217 open instances, at least 8 can be proven optimal without relying on a solver for the proof (\textbf{Table~\ref{tab:solver-free-instances}}). Note that these instances are not simple combinatorial instances, and some are large-scale generic MIPs.
Therefore, one may wonder whether LLMs alone can solve MIPs without solver access.

\begin{table}[ht]
  \centering
  \caption{Eight open \texttt{MIPLIB}
    instances with optimality certificates that do not rely on
    an external solver.}
  \footnotesize
  \renewcommand{\arraystretch}{.94}
  \captionsetup{width=.91\linewidth,skip=4pt}
  \label{tab:solver-free-instances}
  \begin{tabular*}{.98\linewidth}{@{\extracolsep{\fill}}lrrrrrrr@{}}
    \toprule
    Instance & Rows & Columns & Nonzeros & Binary & General int.
      & Cont. & Indicator \\
    \midrule
    \inst{cvrpp-n16k8vrpi}
      & 31,615 & 18,294 & 60,544 & 2,832 & 531 & 14,931 & 13,322 \\
    \inst{dfn-bwin-DBE}
      & 235 & 3,285 & 9,855 & 2,475 & 0 & 810 & 0 \\
    \inst{neos-2629914-sudost}
      & 51,872 & 496 & 209,792 & 256 & 0 & 240 & 0 \\
    \inst{neos-3009394-lami}
      & 2,028 & 2,757 & 6,760 & 2,704 & 52 & 1 & 0 \\
    \inst{neos-3208254-reiu}
      & 9,464,908 & 3,154,336 & 25,080,200 & 3,743 & 0 & 3,150,593 & 0 \\
    \inst{neos-4355351-swalm}
      & 21,609 & 21,065 & 371,467 & 10,530 & 0 & 10,535 & 0 \\
    \inst{neos-4647032-veleka}
      & 8,382 & 18,200 & 3,970,190 & 12,600 & 0 & 5,600 & 0 \\
    \inst{pizza27i}
      & 359,835 & 311,316 & 642,512 & 115,831 & 195,485 & 0 & 97,742 \\
    \bottomrule
  \end{tabular*}
\end{table}

\textbf{Experimental setup.} To test LLMs' potential for solving MIPs directly, we evaluate them on 20 instances from the \texttt{MIPLIB} collection, where each instance is anonymized by permuting rows and columns and renaming variable and constraint names. LLMs can access basic linear algebra libraries such as \texttt{NumPy} and \texttt{SciPy}. Each instance has a time limit of 7200 seconds. The benchmark uses \texttt{COPT 8.0.3} runs with the same time limit.

\textbf{Figure~\ref{fig:solver-replacement}} summarizes the comparison. Without a solver, the LLM can still solve nine instances. Of the 11 that reach the time limit, six finish with a normalized gap below $1\%$. \texttt{COPT} runs report optimality on all 20 instances. This result suggests that for certain MIPs, LLMs alone suffice to solve the MIP instances. However, solvers remain a powerful tool for LLMs to solve MIPs, and we do not recommend prohibiting their use.

\vspace{10pt}

\begin{figure}[ht]
  \centering
  \definecolor{replacementllm}{HTML}{366680}
  \definecolor{replacementcopt}{HTML}{A45B34}
  \begin{tikzpicture}[
    x=.00165cm, y=.14cm, font=\footnotesize,
    llm curve/.style={draw=replacementllm,line width=1.1pt},
    copt curve/.style={draw=replacementcopt,line width=1.1pt,
      dash pattern=on 3pt off 1.6pt},
  ]
    \foreach \level in {0,5,10,15,20} {
      \draw[black!16] (0,\level) -- (7200,\level);
      \node[anchor=east] at (0,\level) {\level};
    }
    \draw[black!65] (0,20) -- (0,0) -- (7200,0);
    \foreach \elapsed in {0,1800,3600,5400,7200} {
      \draw[black!65] (\elapsed,0) -- +(0,-.55);
      \node[anchor=north] at (\elapsed,-.55)
        {\pgfmathprintnumber[1000 sep={,}]{\elapsed}};
    }
    \node[rotate=90] at (-620,10) {Instances solved};
    \node at (3600,-3.9) {Elapsed time per instance (s)};

    \draw[llm curve] (0,0)
      \foreach \elapsed [count=\solved] in
        {150.7,614.6,689,725.5,1165,2540.2,2732.5,2950.2,4792}
        {-| (\elapsed,\solved)} -- (7200,9);
    \draw[copt curve] (0,0)
      \foreach \elapsed [count=\solved] in
        {76,81,109,125,130,699,713,722,922,997,1010,1197,1467,
          1655,1922,2063,2497,3757,3886,5802}
        {-| (\elapsed,\solved)} -- (7200,20);
    \node[anchor=west,text=replacementllm] at (7230,9) {9/20};
    \node[anchor=west,text=replacementcopt] at (7230,20) {20/20};

    \draw[llm curve] (0,23.3) -- (320,23.3)
      node[anchor=west,inner xsep=0pt,xshift=4pt] (llm-legend)
        {LLM without solver (mean time $4{,}777.99$ s)};
    \draw[copt curve] ([xshift=20pt]llm-legend.east) -- ++(320,0)
      node[anchor=west,inner xsep=0pt,xshift=4pt]
        {\texttt{COPT} (mean time $1{,}491.50$ s)};
  \end{tikzpicture}
  \caption{Performance profile on 20
    anonymized \texttt{MIPLIB} instances.}
  \label{fig:solver-replacement}
\end{figure}
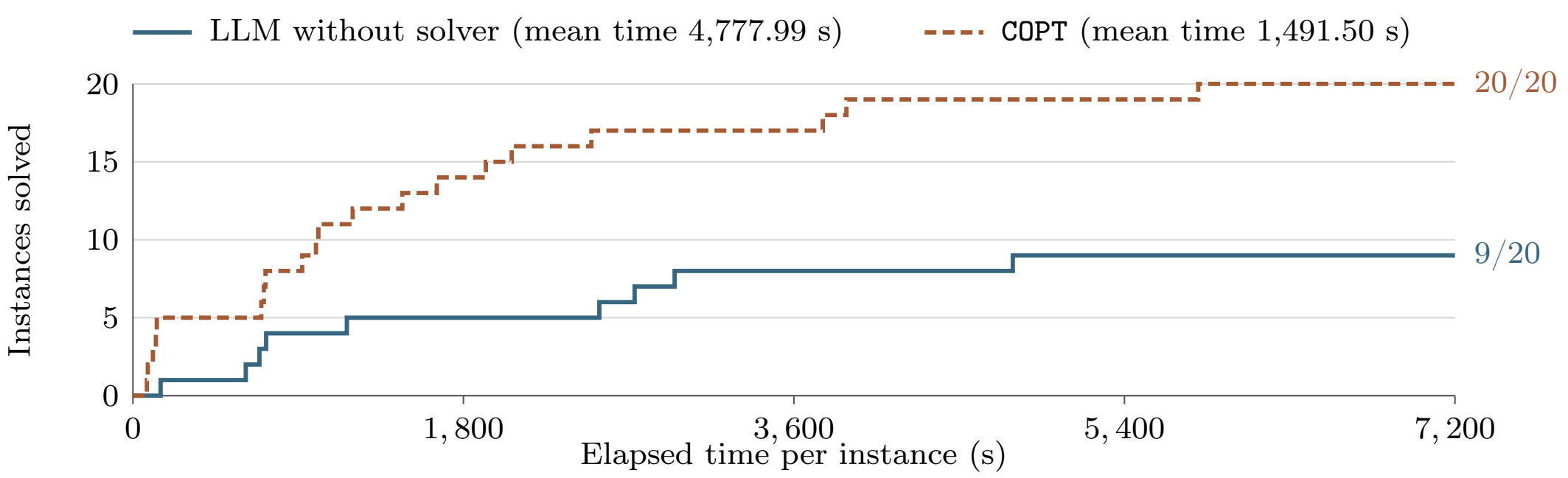

\section{Discussions and conclusions} \label{sec:discussion}

We present {\llm}, the first systematic benchmark of LLMs' ability to exploit structures in MIP instances. Our findings are encouraging: LLMs can already handle complex MIPs and effectively identify useful instance-specific structures. However, we must also point out several limitations of the current LLM-based approaches. \vspace{1pt}

\textbf{Distribution of open instances.} As noted in
\cite{rothberg2026llmsmip}, open instances may represent only part of the range
of MIPs solved in practice. In our experiments, many dual-bound improvements
rely on specialized methods that may apply only to a narrow family of
instances. This limited scope can make the cost of implementing such methods
in a commercial solver hard to justify. \vspace{5pt}

\textbf{Hardness of dual certificates.} Unlike primal bound improvements, which are easy to verify, most dual-side improvements cannot be trivially justified. Unlike systematic frameworks such as branch-and-cut and other relaxation-based approaches, obtaining a certificate by mathematical reasoning or theorem proof is uncommon and still requires massive human intervention. We envision that future MIP solving may finally integrate these non-numerical certificates, but doing it reliably still requires more effort. Some recent attempts have made progress in this direction \cite{robbins2026flare, bai2026rpx}.\vspace{5pt}

\textbf{Privacy.} A thorough investigation of a MIP instance can require
sharing its data or formulation with an LLM. When the data are sensitive, or
the formulation is proprietary, sharing them with a commercial LLM provider
can raise privacy and confidentiality concerns.

\vspace{10pt}

Finally, we discuss the possible implications of LLMs for future solver development. For users, LLMs are likely to significantly lower the barrier for non-experts to solve complicated real-world MIPs. In some cases, an LLM could take on work that would otherwise require a
team of algorithm engineers. Users can easily customize an open-source solver to their needs, or even write a solver from scratch.

From a solver-development perspective, solver development is likely to be more effective. Solver development is similar to the skill distillation process used in the experiments, where methods discovered while solving individual instances can be refined and reused on related problems. A key challenge is to coordinate these methods within a solver while keeping it reliable and easy to use. For example, solvers may develop more callback interfaces compatible with complex structure exploitation.
\vspace{10pt}

Can LLMs help solve MIPs? Our answer is ``Yes, and likely a lot''.

\par
\endgroup

\FloatBarrier
\newpage
\bibliographystyle{plainnat}
\bibliography{references}

\clearpage
\appendix
\section{Confidence score after initial screening}
\label{app:screening-scores}
\noindent
\begin{minipage}{\linewidth}
\centering
\captionsetup{type=table,position=top,skip=5pt}
\caption[\texttt{MIPLIB} 2017 open-instance ranking]{Ranking of 217 open
\texttt{MIPLIB} 2017 instances. Green marks scores $\ge 60$,
yellow 45--59, and gray $<45$.}
\label{tab:miplib-open-ranking}
\begingroup
\definecolor{miprankinggreen}{HTML}{E5F4E8}
\definecolor{miprankingyellow}{HTML}{FFF4D1}
\definecolor{miprankinggray}{HTML}{F0F1F3}
\newcommand*{\miprankingheading}[1]{{\fontsize{7}{8}\selectfont\bfseries #1}}
\newcommand*{\miprankingmetricgap}{3.5pt}
\sffamily\fontsize{8}{7.05}\selectfont
\edef\miprankingpanelwidth{\the\dimexpr(\linewidth-18pt)/3\relax}
\sbox0{\miprankingheading{Primal}\hspace{\miprankingmetricgap}\miprankingheading{Dual}\hspace{\miprankingmetricgap}}
\edef\miprankingbandwidth{\the\dimexpr\miprankingpanelwidth-\wd0\relax}
\let\miprankingname\miprankingentry
\renewcommand*{\miprankingentry}[2]{\begingroup
  \ifnum#2<45\relax
    \def\miprankingshade{miprankinggray}\else\ifnum#2<60\relax
    \def\miprankingshade{miprankingyellow}\else
    \def\miprankingshade{miprankinggreen}\fi\fi
  \rlap{\smash{\kern-1.5pt\color{\miprankingshade}\rule[-1.15pt]{\miprankingbandwidth}{6.65pt}}}\endgroup
  \miprankingname{#1} & #2\kern1.5pt}
\newenvironment{miprankingpanel}{\begin{tabular*}{\miprankingpanelwidth}[t]{@{\hspace{1.5pt}}l@{\extracolsep{\fill}\hspace{1pt}}r@{\extracolsep{0pt}\hspace{\miprankingmetricgap}}r@{\hspace{\miprankingmetricgap}}r@{}}
  \toprule
  \miprankingheading{Instance} & \miprankingheading{Score} &
  \miprankingheading{Primal} & \miprankingheading{Dual} \\
  \midrule
}{\bottomrule
  \end{tabular*}}
\noindent
\begin{miprankingpanel}
  \miprankingentry{neos-2629914-sudost}{99} & 15 & 99 \\
  \miprankingentry{pythago7825}{98} & 5 & 99 \\
  \miprankingentry{a2864-99blp}{92} & 20 & 96 \\
  \miprankingentry{allcolor58}{82} & 92 & 68 \\
  \miprankingentry{genus-sym-g31-8}{77} & 35 & 82 \\
  \miprankingentry{genus-g31-8}{75} & 35 & 80 \\
  \miprankingentry{graph40-20-1rand}{70} & 35 & 76 \\
  \miprankingentry{supportcase39}{67} & 60 & 66 \\
  \miprankingentry{circ10-3}{66} & 40 & 72 \\
  \miprankingentry{ex1010-pi}{65} & 42 & 70 \\
  \miprankingentry{f2000}{64} & 38 & 70 \\
  \miprankingentry{graph40-40-1rand}{64} & 35 & 70 \\
  \miprankingentry{ramos3}{63} & 38 & 69 \\
  \miprankingentry{genus-sym-grafo5708-48}{62} & 35 & 69 \\
  \miprankingentry{genus-g61-25}{61} & 35 & 68 \\
  \miprankingentry{bppc6-06}{60} & 61 & 57 \\
  \miprankingentry{genus-sym-g62-2}{60} & 35 & 67 \\
  \miprankingentry{assign1-10-4}{59} & 55 & 60 \\
  \miprankingentry{d20200}{58} & 42 & 63 \\
  \miprankingentry{comp12-2idx}{58} & 34 & 59 \\
  \miprankingentry{milo-v13-4-3d-4-0}{57} & 57 & 55 \\
  \miprankingentry{fhnw-binpack4-58}{56} & 66 & 42 \\
  \miprankingentry{neos-3603137-hoteo}{55} & 65 & 40 \\
  \miprankingentry{graph40-80-1rand}{54} & 34 & 62 \\
  \miprankingentry{dfn-bwin-DBE}{54} & 56 & 53 \\
  \miprankingentry{fhnw-binschedule2}{54} & 34 & 49 \\
  \miprankingentry{neos-3691541-lonja}{52} & 24 & 49 \\
  \miprankingentry{neos-4360552-sangro}{52} & 24 & 59 \\
  \miprankingentry{neos-5266653-tugela}{52} & 24 & 49 \\
  \miprankingentry{ns1631475}{52} & 24 & 49 \\
  \miprankingentry{scpj4scip}{52} & 24 & 53 \\
  \miprankingentry{graphdraw-mainerd}{51} & 24 & 43 \\
  \miprankingentry{neos-1420790}{51} & 28 & 53 \\
  \miprankingentry{neos-3594536-henty}{51} & 24 & 54 \\
  \miprankingentry{sorrell7}{51} & 44 & 58 \\
  \miprankingentry{xmas10}{50} & 42 & 57 \\
  \miprankingentry{scpk4}{50} & 20 & 49 \\
  \miprankingentry{xmas10-2}{49} & 42 & 56 \\
  \miprankingentry{neos-4533806-waima}{49} & 24 & 49 \\
  \miprankingentry{neos-4555749-wards}{49} & 24 & 49 \\
  \miprankingentry{t1717}{49} & 34 & 53 \\
  \miprankingentry{t1722}{49} & 34 & 53 \\
  \miprankingentry{bab3}{48} & 30 & 55 \\
  \miprankingentry{lr1dr02vc05v8a-t360}{48} & 34 & 43 \\
  \miprankingentry{neos-3355120-tarago}{48} & 29 & 51 \\
  \miprankingentry{fhnw-schedule-paira200}{48} & 76 & 34 \\
  \miprankingentry{core2586-950}{47} & 24 & 43 \\
  \miprankingentry{core4284-1064}{47} & 24 & 43 \\
  \miprankingentry{core4872-1529}{47} & 24 & 43 \\
  \miprankingentry{fhnw-binschedule0}{47} & 30 & 45 \\
  \miprankingentry{pb-gfrd-pnc}{47} & 24 & 48 \\
  \miprankingentry{pb-grow22}{47} & 24 & 48 \\
  \miprankingentry{scpl4}{47} & 20 & 49 \\
  \miprankingentry{tokyometro}{47} & 28 & 57 \\
  \miprankingentry{fhnw-schedule-paira400}{47} & 76 & 34 \\
  \miprankingentry{fhnw-schedule-pairb200}{47} & 76 & 34 \\
  \miprankingentry{neos-3355323-arnon}{46} & 58 & 35 \\
  \miprankingentry{dws012-01}{46} & 34 & 49 \\
  \miprankingentry{fhnw-schedule-paira100}{46} & 34 & 38 \\
  \miprankingentry{ger50-17-ptp-pop-3t}{46} & 38 & 42 \\
  \miprankingentry{ger50-17-trans-dfn-3t}{46} & 34 & 38 \\
  \miprankingentry{ger50-17-trans-pop-3t}{46} & 34 & 38 \\
  \miprankingentry{graphdraw-grafo2}{46} & 20 & 39 \\
  \miprankingentry{graphdraw-opmanager}{46} & 24 & 43 \\
  \miprankingentry{liu}{46} & 38 & 42 \\
  \miprankingentry{lr2-22dr3-333vc4v17a-t60}{46} & 30 & 45 \\
  \miprankingentry{milo-v12-6-r1-75-1}{46} & 34 & 44 \\
  \miprankingentry{neos-1420546}{46} & 24 & 49 \\
  \miprankingentry{rmine11}{46} & 24 & 48 \\
  \miprankingentry{sing5}{46} & 20 & 45 \\
  \miprankingentry{fhnw-schedule-pairb400}{46} & 76 & 34 \\
  \miprankingentry{gmut-76-40}{45} & 34 & 43 \\
  \miprankingentry{neos-3009394-lami}{45} & 28 & 47 \\
\end{miprankingpanel}\hfill\begin{miprankingpanel}
  \miprankingentry{neos-3695882-vesdre}{45} & 20 & 45 \\
  \miprankingentry{nj1}{44} & 61 & 33 \\
  \miprankingentry{rmine13}{44} & 20 & 44 \\
  \miprankingentry{sing11}{44} & 20 & 45 \\
  \miprankingentry{bley\_xs1}{44} & 63 & 43 \\
  \miprankingentry{nj2}{43} & 60 & 33 \\
  \miprankingentry{lr1dr04vc05v17a-t360}{43} & 34 & 43 \\
  \miprankingentry{ns1456591}{43} & 28 & 47 \\
  \miprankingentry{sct1}{43} & 24 & 44 \\
  \miprankingentry{z26}{43} & 30 & 49 \\
  \miprankingentry{bley\_xs1noM}{43} & 63 & 43 \\
  \miprankingentry{nj3}{42} & 59 & 32 \\
  \miprankingentry{ds-big}{42} & 30 & 49 \\
  \miprankingentry{eva1aprime5x5opt}{42} & 24 & 43 \\
  \miprankingentry{gmut-76-50}{42} & 34 & 43 \\
  \miprankingentry{neos-4531126-vouga}{42} & 25 & 47 \\
  \miprankingentry{scpm1}{42} & 20 & 49 \\
  \miprankingentry{scpn2}{42} & 16 & 45 \\
  \miprankingentry{supportcase22}{42} & 20 & 49 \\
  \miprankingentry{nag}{41} & 28 & 47 \\
  \miprankingentry{neos-2991472-kalu}{41} & 24 & 43 \\
  \miprankingentry{neos-5221106-oparau}{41} & 24 & 43 \\
  \miprankingentry{ns1856153}{41} & 24 & 43 \\
  \miprankingentry{rmine15}{41} & 20 & 44 \\
  \miprankingentry{hgms30}{40} & 30 & 39 \\
  \miprankingentry{neos-3682128-sandon}{40} & 24 & 53 \\
  \miprankingentry{ns1828997}{40} & 24 & 48 \\
  \miprankingentry{sct5}{40} & 24 & 44 \\
  \miprankingentry{set3-09}{40} & 28 & 42 \\
  \miprankingentry{set3-16}{40} & 28 & 42 \\
  \miprankingentry{shipsched}{40} & 34 & 38 \\
  \miprankingentry{ns1905797}{39} & 61 & 27 \\
  \miprankingentry{eva1aprime6x6opt}{39} & 24 & 43 \\
  \miprankingentry{ivu06}{39} & 30 & 49 \\
  \miprankingentry{neos-5045105-creuse}{39} & 28 & 42 \\
  \miprankingentry{rocII-10-11}{39} & 24 & 43 \\
  \miprankingentry{rocII-8-11}{39} & 24 & 43 \\
  \miprankingentry{lr1dr12vc10v70b-t360}{38} & 30 & 39 \\
  \miprankingentry{ns2122698}{38} & 20 & 45 \\
  \miprankingentry{ns930473}{38} & 24 & 43 \\
  \miprankingentry{supportcase23}{38} & 24 & 44 \\
  \miprankingentry{dano3mip}{37} & 24 & 38 \\
  \miprankingentry{dc1l}{37} & 34 & 38 \\
  \miprankingentry{neos-848198}{37} & 24 & 38 \\
  \miprankingentry{nsr8k}{37} & 34 & 38 \\
  \miprankingentry{r4l4-02-tree-bounds-50}{37} & 34 & 38 \\
  \miprankingentry{neos-4562542-watut}{36} & 20 & 39 \\
  \miprankingentry{cmflsp40-36-2-10}{35} & 24 & 38 \\
  \miprankingentry{neos-1423785}{35} & 24 & 38 \\
  \miprankingentry{neos-3634244-kauru}{35} & 20 & 39 \\
  \miprankingentry{ns2124243}{35} & 20 & 45 \\
  \miprankingentry{polygonpack4-10}{35} & 39 & 35 \\
  \miprankingentry{polygonpack4-7}{35} & 39 & 35 \\
  \miprankingentry{ns1690781}{34} & 24 & 43 \\
  \miprankingentry{ns1904248}{34} & 20 & 39 \\
  \miprankingentry{siena1}{34} & 39 & 35 \\
  \miprankingentry{n3705}{34} & 35 & 44 \\
  \miprankingentry{n3707}{34} & 35 & 44 \\
  \miprankingentry{dws012-02}{33} & 34 & 43 \\
  \miprankingentry{dws012-03}{33} & 34 & 43 \\
  \miprankingentry{ivu06-big}{33} & 26 & 45 \\
  \miprankingentry{ivu59}{33} & 26 & 45 \\
  \miprankingentry{neos-2978205-isar}{33} & 20 & 34 \\
  \miprankingentry{neos-4292145-piako}{33} & 24 & 38 \\
  \miprankingentry{neos-4647032-veleka}{33} & 24 & 43 \\
  \miprankingentry{seqsolve2short4288}{33} & 20 & 39 \\
  \miprankingentry{seqsolve3\ldots{}excess384}{33} & 20 & 39 \\ \miprankingentry{n3709}{33} & 35 & 44 \\
  \miprankingentry{n370b}{33} & 35 & 44 \\
  \miprankingentry{cdc7-4-3-2}{32} & 96 & 28 \\
  \miprankingentry{cdma}{32} & 33 & 39 \\
  \miprankingentry{cmflsp60-36-2-6}{32} & 24 & 38 \\
  & & & \\
\end{miprankingpanel}\hfill\begin{miprankingpanel}
  \miprankingentry{hgms62}{32} & 30 & 39 \\
  \miprankingentry{neos-4290317-perth}{32} & 29 & 36 \\
  \miprankingentry{neos-5189128-totara}{32} & 20 & 39 \\
  \miprankingentry{rmine21}{32} & 16 & 40 \\
  \miprankingentry{zeil}{32} & 34 & 38 \\
  \miprankingentry{neos-4355351-swalm}{31} & 24 & 38 \\
  \miprankingentry{seqsolve1}{31} & 20 & 39 \\
  \miprankingentry{stockholm}{31} & 24 & 38 \\
  \miprankingentry{supportcase41}{31} & 24 & 38 \\
  \miprankingentry{nb10tb}{30} & 25 & 37 \\
  \miprankingentry{neos-4358725-tagus}{30} & 29 & 35 \\
  \miprankingentry{polygonpack3-15}{30} & 35 & 31 \\
  \miprankingentry{polygonpack4-15}{30} & 35 & 31 \\
  \miprankingentry{polygonpack5-15}{30} & 35 & 31 \\
  \miprankingentry{cvrpp-n16k8vrpi}{29} & 29 & 35 \\
  \miprankingentry{neos-4545615-waita}{29} & 20 & 39 \\
  \miprankingentry{ns1849932}{28} & 20 & 39 \\
  \miprankingentry{shs1014}{28} & 30 & 39 \\
  \miprankingentry{shs1042}{28} & 30 & 39 \\
  \miprankingentry{neos-4230265-orari}{27} & 29 & 35 \\
  \miprankingentry{rmine25}{27} & 16 & 40 \\
  \miprankingentry{supportcase31}{27} & 25 & 37 \\
  \miprankingentry{n3700}{27} & 39 & 41 \\
  \miprankingentry{neos-5044663-wairoa}{26} & 24 & 38 \\
  \miprankingentry{sing17}{26} & 20 & 39 \\
  \miprankingentry{kottenpark09}{25} & 16 & 45 \\
  \miprankingentry{momentum3}{25} & 29 & 35 \\
  \miprankingentry{neos-2974461-ibar}{25} & 25 & 41 \\
  \miprankingentry{neos-4232544-orira}{25} & 25 & 31 \\
  \miprankingentry{neos-5041822-cockle}{25} & 20 & 39 \\
  \miprankingentry{gsvm2rl11}{25} & 37 & 25 \\
  \miprankingentry{supportcase30}{25} & 22 & 45 \\
  \miprankingentry{gasprod1-1}{24} & 29 & 30 \\
  \miprankingentry{gsvm2rl12}{24} & 37 & 25 \\
  \miprankingentry{fhnw-binschedule1}{23} & 26 & 35 \\
  \miprankingentry{neos-4535459-waipa}{23} & 20 & 39 \\
  \miprankingentry{rwth-timetable}{23} & 20 & 39 \\
  \miprankingentry{uccase10}{23} & 20 & 34 \\
  \miprankingentry{neos-4264598-oueme}{22} & 29 & 30 \\
  \miprankingentry{pizza27i}{22} & 25 & 41 \\
  \miprankingentry{pizza78i}{22} & 25 & 41 \\
  \miprankingentry{fastxgemm-n3r22s4t6}{22} & 20 & 35 \\
  \miprankingentry{gsvm2rl9}{22} & 35 & 25 \\
  \miprankingentry{cvrpb-n45k5vrpi}{21} & 25 & 31 \\
  \miprankingentry{minutedispatchstrategy}{21} & 29 & 30 \\
  \miprankingentry{snp-04-052-052}{21} & 25 & 26 \\
  \miprankingentry{supportcase35}{21} & 29 & 30 \\
  \miprankingentry{fastxgemm-n3r23s5t6}{21} & 20 & 35 \\
  \miprankingentry{proteindesign\ldots{}p19}{20} & 21 & 37 \\ \miprankingentry{tpl-tub-ws1617}{20} & 16 & 35 \\
  \miprankingentry{neos-3754224-navua}{19} & 25 & 31 \\
  \miprankingentry{neos-5151569-mologa}{19} & 25 & 31 \\
  \miprankingentry{fastxgemm-n3r21s3t6}{18} & 20 & 35 \\
  \miprankingentry{bmocbd3}{16} & 25 & 31 \\
  \miprankingentry{neos-3615091-sutlej}{16} & 25 & 31 \\
  \miprankingentry{snp-06-004-052}{16} & 25 & 26 \\
  \miprankingentry{snp-10-004-052}{16} & 25 & 26 \\
  \miprankingentry{snp-10-052-052}{16} & 25 & 26 \\
  \miprankingentry{s82}{15} & 31 & 27 \\
  \miprankingentry{cvrpa-n64k9vrpi}{13} & 25 & 31 \\
  \miprankingentry{mining}{13} & 25 & 26 \\
  \miprankingentry{ns1679495}{13} & 25 & 31 \\
  \miprankingentry{neos-3208254-reiu}{11} & 16 & 35 \\
  \miprankingentry{datt256}{10} & 25 & 12 \\
  \miprankingentry{neos-4285819-pedja}{9} & 21 & 27 \\
  \miprankingentry{bharat}{8} & 21 & 22 \\
  \miprankingentry{dlr1}{8} & 16 & 30 \\
  \miprankingentry{dlr2}{8} & 17 & 18 \\
  \miprankingentry{rfds-4-days}{8} & 31 & 22 \\
  \miprankingentry{splan1}{8} & 21 & 27 \\
  \miprankingentry{supportcase38}{8} & 21 & 22 \\
  \miprankingentry{nucorsav}{8} & 24 & 9 \\
  & & & \\
\end{miprankingpanel}\par
\endgroup
\end{minipage}
\par

\clearpage
\section{Instance bounds and improvements}
\label{app:instance-statistics}
\begingroup
\small
\setlength{\tabcolsep}{3.2pt}
\renewcommand{\arraystretch}{.97}
\setlength{\LTcapwidth}{\linewidth}
\setlength{\LTpre}{3pt}
\setlength{\LTpost}{4pt}

\definecolor{mipstatsprimal}{HTML}{FBE5E5}
\definecolor{mipstatsdual}{HTML}{E4EFF9}
\definecolor{mipstatsoptimal}{HTML}{E5F3E8}
\newcommand*{\mipstatskey}[2]{\colorbox{#1}{\strut #2}}
\newcommand*{\mipstatsgreen}{\cellcolor{mipstatsoptimal}}
\newcommand*{\mipstatsred}{\cellcolor{mipstatsprimal}}
\newcommand*{\mipstatsblue}{\cellcolor{mipstatsdual}}

\noindent
\mipstatskey{mipstatsprimal}{Primal improvement}\quad
\mipstatskey{mipstatsdual}{Dual improvement}\quad
\mipstatskey{mipstatsoptimal}{Reported optimal}
\par\smallskip
\par

\newcommand*{\mipstatshead}{Instance &
  \shortstack[r]{Best-known primal\\\texttt{MIPLIB} v36} &
  \shortstack[r]{\texttt{COPT} dual\\10-hour baseline} &
  \shortstack[r]{\texttt{LLM4MIP}\\primal} &
  \shortstack[r]{\texttt{LLM4MIP}\\dual}\\
}

\begin{longtable}{l*{4}{r}}
\caption{Primal and dual bounds for all 217 \texttt{MIPLIB} instances.}
\label{tab:miplib-instance-statistics}\\
\toprule
\mipstatshead
\midrule
\endfirsthead

\multicolumn{5}{@{}l@{}}{\footnotesize \tablename~\thetable\ (continued)}\\[3pt]
\toprule
\mipstatshead
\midrule
\endhead

\midrule
\multicolumn{5}{@{}r@{}}{\footnotesize Continued on next page}\\
\endfoot

\bottomrule
\endlastfoot
\mipstatsgreen \miprankingentry{a2864-99blp} & $-257$ & $-282.960474$ & $-257$ & \mipstatsblue $-257$ \\
\mipstatsgreen \miprankingentry{allcolor58} & $258$ & $42$ & \mipstatsred $42$ & $42$ \\
\mipstatsgreen \miprankingentry{assign1-10-4} & $422$ & $377.444386$ & $422$ & \mipstatsblue $422$ \\
\miprankingentry{bab3} & $-656214.9542$ & $-657093.783$ & $-656214.9542$ & \textemdash \\
\miprankingentry{bharat} & $4039608.5067$ & $3783727.85$ & \mipstatsred $4012280.0882$ & \mipstatsblue $3850713.4697$ \\
\miprankingentry{bley\_xs1} & $3855895.86$ & $3696540.8$ & $3855895.86$ & $3532608.7406$ \\
\miprankingentry{bley\_xs1noM} & $3873690.77$ & $3713869.16$ & \mipstatsred $3855895.86$ & $3344924.0428$ \\
\miprankingentry{bmocbd3} & $-373568021.47$ & $-373603745$ & \mipstatsred $-373587239.24$ & \mipstatsblue $-373603474.46$ \\
\mipstatsgreen \miprankingentry{bppc6-06} & $208$ & $201$ & $208$ & \mipstatsblue $208$ \\
\miprankingentry{cdc7-4-3-2} & $-307$ & $-381$ & $-307$ & $-381$ \\
\miprankingentry{cdma} & $-2.48\times 10^{16}$ & \textemdash & $-2.47784784\times 10^{16}$ & $-2.700669263\times 10^{16}$ \\
\mipstatsgreen \miprankingentry{circ10-3} & $256$ & $141$ & \mipstatsred $242$ & \mipstatsblue $242$ \\
\miprankingentry{cmflsp40-36-2-10} & $66452234.495$ & $64134678.7$ & $66452236.454$ & \mipstatsblue $65928415.307$ \\
\miprankingentry{cmflsp60-36-2-6} & $73891245.385$ & $72075510.1$ & $73891244.222$ & \mipstatsblue $73054818.854$ \\
\miprankingentry{comp12-2idx} & $277$ & $180.108268$ & $277$ & $165$ \\
\miprankingentry{core2586-950} & $944$ & $937.490655$ & $944$ & $937$ \\
\miprankingentry{core4284-1064} & $1062$ & $1055.87389$ & $1062$ & $1055$ \\
\miprankingentry{core4872-1529} & $1530$ & $1514.10482$ & $1530$ & $932$ \\
\mipstatsgreen \miprankingentry{cvrpa-n64k9vrpi} & $1617$ & $0$ & \mipstatsred $1401$ & \mipstatsblue $1401$ \\
\mipstatsgreen \miprankingentry{cvrpb-n45k5vrpi} & $775$ & $0$ & \mipstatsred $751$ & \mipstatsblue $751$ \\
\mipstatsgreen \miprankingentry{cvrpp-n16k8vrpi} & $450$ & $282.375839$ & $450$ & \mipstatsblue $450$ \\
\miprankingentry{d20200} & $12237$ & $12232.0097$ & $12237$ & $12232$ \\
\miprankingentry{dano3mip} & $664.02083333$ & $578.852969$ & $664.02083333$ & $578.21718558$ \\
\miprankingentry{datt256} & \textemdash & $256$ & \textemdash & \textemdash \\
\miprankingentry{dc1l} & $1758647.6801$ & $1751049.95$ & $1758647.6801$ & $1748926.0308$ \\
\mipstatsgreen \miprankingentry{dfn-bwin-DBE} & $73623.787887$ & $56406.0914$ & $73623.79$ & \mipstatsblue $73623.79$ \\
\miprankingentry{dlr1} & $2708064.137$ & \textemdash & $2708064.1319$ & $2707908.7434$ \\
\miprankingentry{dlr2} & $18311711.897$ & \textemdash & \mipstatsred $18265839.45$ & $11601762.346$ \\
\miprankingentry{ds-big} & $194.58674798$ & $87.2719112$ & $194.49841415$ & \mipstatsblue $178.14457882$ \\
\miprankingentry{dws012-01} & $82030.530143$ & $47506.9179$ & $82030.552222$ & $44641.723169$ \\
\miprankingentry{dws012-02} & $119893.31071$ & $58131.4504$ & \mipstatsred $119443.31072$ & $23459.584922$ \\
\miprankingentry{dws012-03} & $140473.57649$ & $63290.0397$ & $140473.62939$ & \mipstatsblue $114206.30519$ \\
\miprankingentry{eva1aprime5x5opt} & $-17.948075604$ & $-52.0397187$ & $-17.948075604$ & $-458.91907206$ \\
\miprankingentry{eva1aprime6x6opt} & $-18.10146979$ & $-140.895641$ & $-18.101002981$ & \mipstatsblue $-96.703805397$ \\
\miprankingentry{ex1010-pi} & $233$ & $224.057416$ & $233$ & $224$ \\
\miprankingentry{f2000} & $1810$ & $1342.95881$ & $1810$ & $1332$ \\
\miprankingentry{fastxgemm-n3r21s3t6} & $4110.9980012$ & $84$ & \mipstatsred $2147$ & \mipstatsblue $93$ \\
\miprankingentry{fastxgemm-n3r22s4t6} & $3101.9981984$ & $84$ & \mipstatsred $1150$ & \mipstatsblue $90$ \\
\miprankingentry{fastxgemm-n3r23s5t6} & $3089.997998$ & $84$ & $153$ & \mipstatsblue $90$ \\
\miprankingentry{fhnw-binpack4-58} & \textemdash & $0$ & \textemdash & \textemdash \\
\mipstatsgreen \miprankingentry{fhnw-binschedule0} & $16088$ & $15642$ & \mipstatsred $15958$ & \mipstatsblue $15958$ \\
\mipstatsgreen \miprankingentry{fhnw-binschedule1} & $55198$ & $52153$ & \mipstatsred $55158$ & \mipstatsblue $55158$ \\
\mipstatsgreen \miprankingentry{fhnw-binschedule2} & $2428$ & $2273.5$ & $2428$ & \mipstatsblue $2428$ \\
\mipstatsgreen \miprankingentry{fhnw-schedule-paira100} & $-15.113122668$ & $-23.8016842$ & $-15.113116513$ & \mipstatsblue $-15.113116513$ \\
\mipstatsgreen \miprankingentry{fhnw-schedule-paira200} & $-19.229856822$ & $-46.9081453$ & \mipstatsred $-19.369961205$ & \mipstatsblue $-19.369961205$ \\
\mipstatsgreen \miprankingentry{fhnw-schedule-paira400} & $-35.546804653$ & $-95.8844179$ & \mipstatsred $-36.276966742$ & \mipstatsblue $-36.276966742$ \\
\mipstatsgreen \miprankingentry{fhnw-schedule-pairb200} & $-19.240943843$ & $-51.5610835$ & \mipstatsred $-19.369961205$ & \mipstatsblue $-19.369961205$ \\
\mipstatsgreen \miprankingentry{fhnw-schedule-pairb400} & $-35.4584$ & $-98.8384308$ & \mipstatsred $-36.276966742$ & \mipstatsblue $-36.276966742$ \\
\mipstatsgreen \miprankingentry{gasprod1-1} & $291.17532062$ & $253.271446$ & $291.52653003$ & \mipstatsblue $291.52653003$ \\
\mipstatsgreen \miprankingentry{genus-g31-8} & $-23$ & $-27$ & $-23$ & \mipstatsblue $-23$ \\
\miprankingentry{genus-g61-25} & $-40$ & $-58$ & $-40$ & \mipstatsblue $-42$ \\
\mipstatsgreen \miprankingentry{genus-sym-g31-8} & $-23$ & $-27$ & $-23$ & \mipstatsblue $-23$ \\
\miprankingentry{genus-sym-g62-2} & $-38$ & $-54$ & $-38$ & \mipstatsblue $-40$ \\
\mipstatsgreen \miprankingentry{genus-sym-grafo5708-48} & $-21$ & $-27$ & $-21$ & \mipstatsblue $-21$ \\
\miprankingentry{ger50-17-ptp-pop-3t} & $5224.5144$ & $4123.52638$ & $5224.5144$ & \mipstatsblue $5183.9359259$ \\
\miprankingentry{ger50-17-trans-dfn-3t} & $3969.4334$ & $3831.36744$ & $3969.4334$ & \mipstatsblue $3934.3152$ \\
\mipstatsgreen \miprankingentry{ger50-17-trans-pop-3t} & $4038.441$ & $3884.84515$ & $4038.441$ & \mipstatsblue $4038.441$ \\
\miprankingentry{gmut-76-40} & $-14169477.796$ & $-14170137.1$ & $-14169477.786$ & $-14171086.316$ \\
\miprankingentry{gmut-76-50} & $-14171893.779$ & $-14173526.5$ & $-14171397.828$ & \textemdash \\
\miprankingentry{graph40-20-1rand} & $-15$ & $-134.42884$ & $-15$ & \textemdash \\
\mipstatsgreen \miprankingentry{graph40-40-1rand} & $-9$ & $-68.9559386$ & $-9$ & \mipstatsblue $-9$ \\
\mipstatsgreen \miprankingentry{graph40-80-1rand} & $-7$ & $-495.450445$ & $-7$ & \mipstatsblue $-7$ \\
\miprankingentry{graphdraw-grafo2} & $68613.5$ & $36314.1643$ & \mipstatsred $68590.5$ & \mipstatsblue $53513.5$ \\
\miprankingentry{graphdraw-mainerd} & $38842$ & $31493.2257$ & $38842$ & \mipstatsblue $31882$ \\
\miprankingentry{graphdraw-opmanager} & $99445.5$ & $56741.4457$ & $99445.5$ & \mipstatsblue $58950$ \\
\miprankingentry{gsvm2rl11} & $18121.637644$ & $3169.43538$ & $18539.127162$ & $2096.7554578$ \\
\miprankingentry{gsvm2rl12} & $21.903856359$ & $3.11194162$ & $21.903856359$ & \mipstatsblue $12.099305681$ \\
\miprankingentry{gsvm2rl9} & $7438.1810212$ & $2889.25545$ & $7438.1811678$ & $2289.6595304$ \\
\miprankingentry{hgms30} & $-44338.744803$ & $-46112.9057$ & $-44338.744813$ & \mipstatsblue $-44339.653641$ \\
\miprankingentry{hgms62} & $-44837.245343$ & $-46639.631$ & $-44837.245347$ & \mipstatsblue $-44844.328684$ \\
\miprankingentry{ivu06} & $142.86$ & $135.944621$ & $142.86433158$ & $135.90024232$ \\
\miprankingentry{ivu06-big} & $140.74$ & $135.430283$ & \mipstatsred $140.41233208$ & \mipstatsblue $135.44891253$ \\
\miprankingentry{ivu59} & $927.89475182$ & $884.788371$ & \mipstatsred $911.94105059$ & \mipstatsblue $885$ \\
\miprankingentry{kottenpark09} & $1715$ & $183.5$ & $1715$ & \mipstatsblue $230$ \\
\miprankingentry{liu} & $1083.984782$ & $560$ & \mipstatsred $1082$ & \mipstatsblue $1074$ \\
\miprankingentry{lr1dr02vc05v8a-t360} & $123046.81494$ & $120969.134$ & $123046.81612$ & \textemdash \\
\miprankingentry{lr1dr04vc05v17a-t360} & $252463$ & $247825.753$ & \mipstatsred $251302.15104$ & \mipstatsblue $247852.13666$ \\
\miprankingentry{lr1dr12vc10v70b-t360} & $1000790.9797$ & $984474.155$ & \mipstatsred $999429.32495$ & \mipstatsblue $985808.99169$ \\
\miprankingentry{lr2-22dr3-333vc4v17a-t60} & $-27262.589176$ & $-71836.9293$ & \mipstatsred $-27640.919482$ & \mipstatsblue $-30627.159226$ \\
\miprankingentry{milo-v12-6-r1-75-1} & $1153756.3982$ & $1135852.03$ & $1153756.3982$ & $1125907.0416$ \\
\miprankingentry{milo-v13-4-3d-4-0} & $358152.26219$ & $253283.451$ & $358152.26219$ & $221230.41434$ \\
\miprankingentry{mining} & $-833589149.27$ & $-933493233$ & $-833589149.27$ & \mipstatsblue $-835361895.59$ \\
\mipstatsgreen \miprankingentry{minutedispatchstrategy} & $3109.9034778$ & $3109.78686$ & $3109.9034778$ & \mipstatsblue $3109.9034778$ \\
\miprankingentry{momentum3} & $173202.17622$ & $96179.1231$ & $173202.18557$ & \mipstatsblue $107995.45322$ \\
\miprankingentry{n3700} & $1218975$ & $1065546$ & $1247295$ & $1063103.1471$ \\
\miprankingentry{n3705} & $1212657$ & $1064197.11$ & \mipstatsred $1208356$ & \mipstatsblue $1129725$ \\
\miprankingentry{n3707} & $1186691$ & $1035948.62$ & \mipstatsred $1185209$ & \mipstatsblue $1102426$ \\
\miprankingentry{n3709} & $1205801$ & $1054892.55$ & $1205801$ & \mipstatsblue $1146484$ \\
\miprankingentry{n370b} & $1220708$ & $1074178.91$ & \mipstatsred $1212433$ & \mipstatsblue $1166799$ \\
\miprankingentry{nag} & $945$ & $465$ & \mipstatsred $930$ & \mipstatsblue $550.54054054$ \\
\miprankingentry{nb10tb} & $1.3796452423\times 10^{10}$ & $1.33301713\times 10^{10}$ & \mipstatsred $1.3796451833\times 10^{10}$ & \mipstatsblue $1.3794494287\times 10^{10}$ \\
\miprankingentry{neos-1420546} & $22948.103198$ & $22223.7073$ & \mipstatsred $22936.284978$ & \mipstatsblue $22629.329177$ \\
\miprankingentry{neos-1420790} & $3119.8050577$ & $3003.05309$ & $3119.8060016$ & \mipstatsblue $3067.0672911$ \\
\miprankingentry{neos-1423785} & $21893.26$ & $20327.177$ & $21893.333343$ & \mipstatsblue $20707.764318$ \\
\mipstatsgreen \miprankingentry{neos-2629914-sudost} & $48180$ & $42387.0653$ & $48180$ & \mipstatsblue $48180$ \\
\miprankingentry{neos-2974461-ibar} & $468906174.77$ & $429197921$ & $468906174.77$ & \mipstatsblue $443920829.19$ \\
\mipstatsgreen \miprankingentry{neos-2978205-isar} & $-11.92580869$ & $-11.9872276$ & $-11.925808687$ & \mipstatsblue $-11.925808687$ \\
\mipstatsgreen \miprankingentry{neos-2991472-kalu} & $12$ & $0$ & $12$ & \mipstatsblue $12$ \\
\mipstatsgreen \miprankingentry{neos-3009394-lami} & $5.5$ & $1.92857143$ & $5.5$ & \mipstatsblue $5.5$ \\
\mipstatsgreen \miprankingentry{neos-3208254-reiu} & $-34897.000525$ & $-38606.2553$ & \mipstatsred $-34902$ & \mipstatsblue $-34902$ \\
\miprankingentry{neos-3355120-tarago} & $-11115971.494$ & $-11712060.9$ & $-11110906.689$ & \mipstatsblue $-11123832.677$ \\
\miprankingentry{neos-3355323-arnon} & \textemdash & $0$ & \textemdash & \textemdash \\
\mipstatsgreen \miprankingentry{neos-3594536-henty} & $401223$ & $397821.689$ & \mipstatsred $401092$ & \mipstatsblue $401092$ \\
\miprankingentry{neos-3603137-hoteo} & \textemdash & $0$ & \textemdash & \textemdash \\
\mipstatsgreen \miprankingentry{neos-3615091-sutlej} & $-1048730.9379$ & $-1404604.89$ & $-1048730.8$ & \mipstatsblue $-1048730.8$ \\
\miprankingentry{neos-3634244-kauru} & $1398.0146705$ & $239.044636$ & $1398.0156859$ & \textemdash \\
\mipstatsgreen \miprankingentry{neos-3682128-sandon} & $34666767.474$ & $21531270.9$ & $34666770$ & \mipstatsblue $34666770$ \\
\miprankingentry{neos-3691541-lonja} & $-0.00744145058$ & $-0.019851233$ & $0.0918239856$ & \textemdash \\
\miprankingentry{neos-3695882-vesdre} & $518.35054802$ & $-7522.4774$ & $518.5$ & \mipstatsblue $517.4$ \\
\miprankingentry{neos-3754224-navua} & $157909.06776$ & $56251.4076$ & $166345.74719$ & \mipstatsblue $56490.877344$ \\
\miprankingentry{neos-4230265-orari} & $73755$ & $21013.257$ & \mipstatsred $72572$ & \mipstatsblue $31528$ \\
\miprankingentry{neos-4232544-orira} & $5538597.4$ & $5070682.61$ & $5547858.8$ & \mipstatsblue $5099038.4143$ \\
\miprankingentry{neos-4264598-oueme} & $5997393.4777$ & $5388920.94$ & $6207553.6$ & \mipstatsblue $5468401.5318$ \\
\miprankingentry{neos-4285819-pedja} & $88.499336306$ & $63.3949254$ & $88.499336306$ & $62.376611387$ \\
\mipstatsgreen \miprankingentry{neos-4290317-perth} & $3017215.8714$ & $3017032.16$ & \mipstatsred $3017171.1908$ & \mipstatsblue $3017171.1908$ \\
\miprankingentry{neos-4292145-piako} & $19229.998369$ & $15702.1569$ & $27660$ & \mipstatsblue $27570.1$ \\
\mipstatsgreen \miprankingentry{neos-4355351-swalm} & $33.45757454$ & $33.45763$ & $33.45763$ & $33.45763$ \\
\mipstatsgreen \miprankingentry{neos-4358725-tagus} & $4206.1401$ & $2378.09139$ & $32225.700099$ & $32225.700099$ \\
\mipstatsgreen \miprankingentry{neos-4360552-sangro} & $-8$ & $-12.0644768$ & $-8$ & \mipstatsblue $-8$ \\
\miprankingentry{neos-4531126-vouga} & $13216.633968$ & $412238.281$ & $525030.88462$ & \mipstatsblue $524657.38188$ \\
\miprankingentry{neos-4533806-waima} & $4870$ & $4249.09466$ & $4870$ & $3956$ \\
\miprankingentry{neos-4535459-waipa} & $26040609$ & $0$ & \mipstatsred $23331046$ & \mipstatsblue $15474160.819$ \\
\miprankingentry{neos-4545615-waita} & $6311$ & $3753.09803$ & $6311$ & \mipstatsblue $6021$ \\
\mipstatsgreen \miprankingentry{neos-4555749-wards} & $3759$ & $3187.97075$ & $3759$ & \mipstatsblue $3759$ \\
\mipstatsgreen \miprankingentry{neos-4562542-watut} & $5684$ & $4706.80045$ & $5684$ & \mipstatsblue $5684$ \\
\mipstatsgreen \miprankingentry{neos-4647032-veleka} & $27214.4801$ & $27213.9842$ & $27214.4801$ & \mipstatsblue $27214.4801$ \\
\miprankingentry{neos-5041822-cockle} & $5724245.7045$ & $5697141.9$ & $5724245.654$ & \mipstatsblue $5701917.8171$ \\
\miprankingentry{neos-5044663-wairoa} & $4654484.5553$ & $4654134.39$ & $4654500$ & \mipstatsblue $4654300$ \\
\miprankingentry{neos-5045105-creuse} & $20.571247075$ & $20.1974827$ & $20.571410409$ & $19.932229271$ \\
\miprankingentry{neos-5151569-mologa} & $686748326.49$ & $686479073$ & \mipstatsred $686721918.74$ & \mipstatsblue $686532411.38$ \\
\miprankingentry{neos-5189128-totara} & $6459876.1385$ & $501198.807$ & \mipstatsred $2901877.8002$ & \mipstatsblue $2494433.1727$ \\
\miprankingentry{neos-5221106-oparau} & $52.67$ & $0$ & \mipstatsred $52.07$ & \mipstatsblue $42.54$ \\
\miprankingentry{neos-5266653-tugela} & $65505.200575$ & $40074.8618$ & \mipstatsred $65257.57399$ & \textemdash \\
\miprankingentry{neos-848198} & $51797$ & $50021.1801$ & $51797$ & $49843$ \\
\miprankingentry{nj1} & \textemdash & $78.882279$ & \mipstatsred $382.0606208$ & \textemdash \\
\miprankingentry{nj2} & \textemdash & $47.1019816$ & \mipstatsred $382.0606208$ & \textemdash \\
\miprankingentry{nj3} & \textemdash & $26.6759647$ & \mipstatsred $382.0606208$ & \textemdash \\
\mipstatsgreen \miprankingentry{ns1456591} & $988.14054407$ & $459.574181$ & $988.14128344$ & \mipstatsblue $988.14128344$ \\
\mipstatsgreen \miprankingentry{ns1631475} & $11100$ & $1659$ & $11100$ & \mipstatsblue $11100$ \\
\miprankingentry{ns1679495} & $2642035.1987$ & $2556399.7$ & $2642035.1987$ & \mipstatsblue $2626815.6662$ \\
\miprankingentry{ns1690781} & $-928.08534703$ & $-1022.68453$ & \mipstatsred $-928.18882865$ & \mipstatsblue $-932.72130546$ \\
\miprankingentry{ns1828997} & $8$ & $0$ & $8$ & \textemdash \\
\miprankingentry{ns1849932} & $632$ & $546$ & \mipstatsred $624$ & \mipstatsblue $548$ \\
\miprankingentry{ns1856153} & $34.1634618$ & $0$ & \mipstatsred $34.016973126$ & \mipstatsblue $32.697165488$ \\
\miprankingentry{ns1904248} & $37.017626111$ & $0$ & \mipstatsred $35.262535127$ & \mipstatsblue $23.6999997$ \\
\miprankingentry{ns1905797} & \textemdash & $12.38858$ & \mipstatsred $13.98258$ & $12.38858$ \\
\miprankingentry{ns2122698} & $68749.999885$ & $65899.0548$ & $68750$ & $62785.565316$ \\
\mipstatsgreen \miprankingentry{ns2124243} & $77459.991691$ & $75440.7297$ & $77460$ & \mipstatsblue $77460$ \\
\mipstatsgreen \miprankingentry{ns930473} & $821466$ & $1899.50067$ & $821466$ & \mipstatsblue $821466$ \\
\miprankingentry{nsr8k} & $17974887.998$ & $17525904.2$ & $17974887.998$ & $17512125.696$ \\
\mipstatsgreen \miprankingentry{nucorsav} & \textemdash & $29.0461$ & \mipstatsred $73.61064226$ & \mipstatsblue $73.610642246$ \\
\miprankingentry{pb-gfrd-pnc} & $8844$ & $8521.83408$ & $8844$ & \mipstatsblue $8566$ \\
\miprankingentry{pb-grow22} & $-449536$ & $-449587.303$ & $-449536$ & \mipstatsblue $-449546$ \\
\mipstatsgreen \miprankingentry{pizza27i} & $701882$ & $652380.8$ & $701882$ & \mipstatsblue $701882$ \\
\mipstatsgreen \miprankingentry{pizza78i} & $564039$ & $508093$ & $564039$ & \mipstatsblue $564039$ \\
\miprankingentry{polygonpack3-15} & $-55055635.466$ & $-114051303$ & \mipstatsred $-55657792.808$ & \mipstatsblue $-58365098.476$ \\
\miprankingentry{polygonpack4-10} & $-53594508.708$ & $-93676573.4$ & $-53594508.708$ & \mipstatsblue $-76521127.377$ \\
\miprankingentry{polygonpack4-15} & $-63613643.562$ & $-114052544$ & \textemdash & \textemdash \\
\miprankingentry{polygonpack4-7} & $-51837758.962$ & $-75324457.9$ & $-51837707.95$ & \textemdash \\
\miprankingentry{polygonpack5-15} & $-55494687.648$ & $-66029891.1$ & \textemdash & \textemdash \\
\mipstatsgreen \miprankingentry{proteindesign121hz512p19} & $3382$ & $3392$ & $3382$ & $3382$ \\
\miprankingentry{pythago7825} & \textemdash & $0$ & \textemdash & \textemdash \\
\miprankingentry{r4l4-02-tree-bounds-50} & $499132179$ & $478217938$ & \mipstatsred $499061405$ & \mipstatsblue $483755247.89$ \\
\miprankingentry{ramos3} & $186$ & $146.482178$ & $186$ & \mipstatsblue $156$ \\
\miprankingentry{rfds-4-days} & $-1129915.8646$ & $-1157062.2$ & $-1129915.8646$ & $-1146278.6542$ \\
\miprankingentry{rmine11} & $-2508.404144$ & $-2509.7041$ & $-2508.404144$ & $-2521.3532743$ \\
\miprankingentry{rmine13} & $-3495.3706624$ & $-3496.99083$ & $-3495.370661$ & $-3504.6899591$ \\
\miprankingentry{rmine15} & $-5018.819991$ & $-5026.70481$ & \mipstatsred $-5018.823005$ & $-5042.0362334$ \\
\miprankingentry{rmine21} & $-10618.750838$ & $-10674.8943$ & $-10618.750727$ & $-10679.232567$ \\
\miprankingentry{rmine25} & $-15541.669287$ & $-15667.8428$ & $-15541.669277$ & $-15667.918863$ \\
\miprankingentry{rocII-10-11} & $-11.793992743$ & $-11.8947324$ & $-11.794324502$ & \mipstatsblue $-11.849392208$ \\
\miprankingentry{rocII-8-11} & $-8.7398453503$ & $-11.889881$ & $-8.739845027$ & \mipstatsblue $-11.82$ \\
\miprankingentry{rwth-timetable} & $772998$ & $754720.911$ & \mipstatsred $772258$ & \mipstatsblue $755608.65978$ \\
\miprankingentry{s82} & $-33.797057624$ & $-33.8093843$ & $-33.797057624$ & \textemdash \\
\miprankingentry{scpj4scip} & $128$ & $108.947237$ & $128$ & $107$ \\
\miprankingentry{scpk4} & $318$ & $287.057116$ & $318$ & $283$ \\
\miprankingentry{scpl4} & $259$ & $214.35471$ & $259$ & $214$ \\
\miprankingentry{scpm1} & $540$ & $415.555435$ & $537$ & \mipstatsblue $416$ \\
\miprankingentry{scpn2} & $485$ & $350.409358$ & $485$ & \mipstatsblue $352$ \\
\miprankingentry{sct1} & $-187.52765084$ & $-202.560756$ & $-187.52765084$ & \mipstatsblue $-199.81870652$ \\
\miprankingentry{sct5} & $-228.12919207$ & $-228.154186$ & $-228.12919207$ & \mipstatsblue $-228.14788057$ \\
\miprankingentry{seqsolve1} & $4277$ & $4268.52039$ & $4277$ & \mipstatsblue $4275$ \\
\miprankingentry{seqsolve2short4288} & $84$ & $33.1143818$ & $83.999999835$ & $17$ \\
\miprankingentry{seqsolve3short4288excess384} & $-40015$ & $-42615.0595$ & $-40035$ & \mipstatsblue $-41168$ \\
\miprankingentry{set3-09} & $176497.15$ & $134902.999$ & \textemdash & \textemdash \\
\miprankingentry{set3-16} & $134040.40975$ & $104243.075$ & $134040.40975$ & $101848.66004$ \\
\miprankingentry{shipsched} & $111919.83373$ & $67029.7019$ & $111919.83373$ & \mipstatsblue $85777.682467$ \\
\miprankingentry{shs1014} & $22671.197427$ & $23116.6794$ & $22671.197442$ & $22555.877512$ \\
\miprankingentry{shs1042} & $11070.70158$ & $11218.5138$ & $11070.701662$ & $10905.739776$ \\
\miprankingentry{siena1} & $10341586.2$ & $10214130.1$ & $10449008.9$ & \mipstatsblue $10219502.195$ \\
\miprankingentry{sing11} & $19108639.439$ & $18776564.1$ & $19108639.46$ & \mipstatsblue $18853382.032$ \\
\miprankingentry{sing17} & $36161699.379$ & $35599182.4$ & $36161699.378$ & $35556219.48$ \\
\miprankingentry{sing5} & $18778435.608$ & $18605207.3$ & $18778435.608$ & $18462676.779$ \\
\miprankingentry{snp-04-052-052} & $881796490.75$ & $843534432$ & \mipstatsred $881740229.2$ & \mipstatsblue $880540652.11$ \\
\miprankingentry{snp-06-004-052} & $1869531794.3$ & $1869529780$ & $1869531794.2$ & \mipstatsblue $1869530183.6$ \\
\miprankingentry{snp-10-004-052} & $5906642759.1$ & $5906557160$ & \mipstatsred $5906642753.4$ & \mipstatsblue $5906641320.2$ \\
\miprankingentry{snp-10-052-052} & $6059860690.3$ & $6011665050$ & \mipstatsred $6059234264.5$ & \mipstatsblue $6055563163.9$ \\
\miprankingentry{sorrell7} & $-198$ & $-223.02692$ & $-198$ & \mipstatsblue $-210$ \\
\miprankingentry{splan1} & $-194377.63065$ & $-201881.851$ & \mipstatsred $-194385.95334$ & \mipstatsblue $-201513.1968$ \\
\miprankingentry{stockholm} & $120.99999502$ & $93.5364881$ & $122$ & \mipstatsblue $106$ \\
\miprankingentry{supportcase22} & $110$ & $0.88353308$ & $110$ & \mipstatsblue $6.4$ \\
\miprankingentry{supportcase23} & $-12160.659357$ & $-12633.1586$ & $-12160.659354$ & \mipstatsblue $-12345.870234$ \\
\miprankingentry{supportcase30} & \textemdash & $0$ & \textemdash & \textemdash \\
\mipstatsgreen \miprankingentry{supportcase31} & $-3720096.597$ & $-3721177.76$ & \mipstatsred $-3720097.298$ & \mipstatsblue $-3720097.298$ \\
\miprankingentry{supportcase35} & $-164208.6431$ & \textemdash & \textemdash & \textemdash \\
\miprankingentry{supportcase38} & $32873.968346$ & $32476.0327$ & \mipstatsred $32564.500163$ & $32476.032696$ \\
\miprankingentry{supportcase39} & $-1085080.9069$ & $-1085222.66$ & \mipstatsred $-1085083.8568$ & $-1085325.2307$ \\
\miprankingentry{supportcase41} & $32.000768321$ & $28$ & $32.000768325$ & \mipstatsblue $31.81114547$ \\
\miprankingentry{t1717} & $158260$ & $137668.992$ & $158260$ & $137323.70815$ \\
\miprankingentry{t1722} & $108953$ & $103989.69$ & $108953$ & \textemdash \\
\miprankingentry{tokyometro} & $8263.1$ & $6817.12644$ & $8263.1$ & \mipstatsblue $7117.5$ \\
\miprankingentry{tpl-tub-ws1617} & $121066$ & $120987.401$ & $121065$ & $120885$ \\
\miprankingentry{uccase10} & $39133.58411$ & $39113.5635$ & $39133.58411$ & \mipstatsblue $39122.060112$ \\
\miprankingentry{xmas10} & $-497$ & $-516.336258$ & $-497$ & \mipstatsblue $-507$ \\
\miprankingentry{xmas10-2} & $-497$ & $-515.926191$ & $-497$ & \mipstatsblue $-515$ \\
\miprankingentry{z26} & $-1198$ & $-1273.9627$ & $-1198$ & $-1646$ \\
\miprankingentry{zeil} & $1081.3585$ & $834.502247$ & $1081.3585$ & $833.83125413$ \\

\end{longtable}

\noindent
\par
\endgroup

\endgroup
\end{document}